\documentclass[12pt]{amsart}
\usepackage{amsmath}
\usepackage{amsfonts}
\usepackage{epsfig}
\usepackage{fullpage,url}
\usepackage[dvips]{hyperref}
 \usepackage{amscd}   

\DeclareFontEncoding{OT2}{}{} 

\newcommand{\pp}{{\mathfrak p}}

\DeclareMathOperator{\Stab}{Stab}

\DeclareMathOperator{\F}{\sf{F}}
\DeclareMathOperator{\GF}{\sf{GF}}
\DeclareMathOperator{\Gaa}{\tilde{\Gamma}}

\DeclareMathOperator{\End}{End}

\DeclareMathOperator{\tor}{tor}
\DeclareMathOperator{\alg}{alg}

\DeclareMathOperator{\Frac}{Frac}

\newcommand{\isom}{\simeq}

\newcommand{\tensor}{\otimes}

\newtheorem{theorem}{Theorem}[section]
\newtheorem{lemma}[theorem]{Lemma}

\theoremstyle{definition}
\newtheorem{definition}[theorem]{Definition}

\newtheorem{example}[theorem]{Example}

\theoremstyle{remark}
\newtheorem{remark}[theorem]{Remark}

\title{The Mordell-Lang Theorem for finitely generated subgroups of a semiabelian variety defined over a finite field}
\author{Dragos Ghioca}
\address{Dragos Ghioca, Department of Mathematics \& Statistics, Hamilton Hall, Room 218, McMaster University, 1280 Main Street West, Hamilton, Ontario L8S 4K1, Canada}

\email{dghioca@math.mcmaster.ca}
\begin{document}

\begin{abstract}
We determine the structure of the intersection of a finitely generated subgroup of a semiabelian variety $G$ defined over a finite field with a closed subvariety $X\subset G$.
\end{abstract}

\maketitle

\section{Introduction}
\label{S:intro}

\footnotetext[1]{2000 AMS Subject Classification: Primary 11G10; Secondary 11G09, 11G25 }

Let $G$ be a semiabelian variety defined over a finite field $\mathbb{F}_q$. Let $K$ be a regular field extension of $\mathbb{F}_q$. Let $F$ be the corresponding Frobenius for $\mathbb{F}_q$. Then $F\in\End(G)$. 

Let $X$ be a subvariety of $G$ defined over $K$ (in this paper, all subvarieties will be closed). In \cite{F-sets} and \cite{revisited}, Moosa and Scanlon discussed the intersection of the $K$-points of $X$ with a finitely generated $\mathbb{Z}[F]$-submodule $\Gamma$ of $G(K)$. They proved that the intersection is a finite union of \emph{$F$-sets} in $\Gamma$ (see Definition~\ref{D:F-sets}). Our goal is to extend their result to the case when $\Gamma$ is a finitely generated subgroup of $G(K)$ (not necessarily invariant under $F$). 

In Section~\ref{S:statement} we will state our main results, which include, besides the Mordell-Lang statement for subgroups of semiabelian varieties defined over finite fields, also a similar Mordell-Lang statement for Drinfeld modules defined over finite fields. The Mordell-Lang Theorem for Drinfeld modules was also studied by the author in \cite{IMRN}. In Section~\ref{S:semiabelian} we will prove our main theorem for semiabelian varieties, while in Section~\ref{S:Drinfeld} we will show how the Mordell-Lang statement for Drinfeld modules defined over finite fields can be deduced from the results in \cite{F-sets}. We will conclude Section~\ref{S:Drinfeld} with two counterexamples for two possible extensions of our statement for Drinfeld modules towards results similar with the ones true for semiabelian varieties.

\section{Statement of our main results}
\label{S:statement}

Everywhere in this paper, $\overline{Y}$ represents the Zariski closure of the set $Y$.

A central notion for the present paper is the notion of a \emph{Frobenius ring}. This notion was first introduced by Moosa and Scanlon (see Definition $2.1$ in \cite{revisited}). We extend their definition to include also rings of finite characteristic.
\begin{definition}
\label{D:Frobenius ring}
Let $R$ be a Dedekind domain with the property that for every nonzero prime ideal $\pp\subset R$, $R/\pp$ is a finite field. We call $R[F]$ a Frobenius ring if the following properties are satisfied:

$(i)$ $R[F]$ is a simple extension of $R$ generated by a distinguished element $F$.

$(ii)$ $R[F]$ is a finite integral extension of $R$.

$(iii)$ $F$ is not a zero divisor in $R[F]$.

$(iv)$ The ideal $F^{\infty}R[F]:=\bigcap_{n\ge 0}F^nR[F]$ is trivial.
\end{definition}

The classical example of a Frobenius ring associated to a semiabelian variety $G$ defined over the finite field $\mathbb{F}_q$ is $\mathbb{Z}[F]$, where $F$ is the corresponding Frobenius for $\mathbb{F}_q$. This Frobenius ring is discussed in \cite{F-sets} and \cite{revisited}. We will show later in this section that $A[F]$ is also a Frobenius ring when $F$ is the Frobenius on $\mathbb{F}_q$ and $\phi:A\rightarrow\mathbb{F}_q[F]$ is a Drinfeld module (in this case, $A$ is a Dedekind domain of finite characteristic).

We define the notion of \emph{groupless $F$-sets} contained in a module over a Frobenius ring.
\begin{definition}
\label{D:groupless}
Let $R[F]$ be a Frobenius ring and let $M$ be an $R[F]$-module. 
For $a\in M$ and $\delta\in\mathbb{N}^{*}$, we denote the $F^{\delta}$-orbit of $a$ by $S(a;\delta):=\{F^{\delta n}a\mid n\in\mathbb{N}\}$. If $a_1,\dots,a_k\in M$ and $\delta_1,\dots,\delta_k\in\mathbb{N}^{*}$, then we denote the sum of the $F^{\delta_i}$-orbits of $a_i$ by $$S(a_1,\dots,a_k;\delta_1,\dots,\delta_k)=\{\sum_{i=1}^kF^{\delta_i n_i}a_i\mid (n_1,\dots,n_k)\in\mathbb{N}^k\}.$$
A set of the form $b+S(a_1,\dots,a_k;\delta_1,\dots,\delta_k)$ with $b,a_1,\dots,a_k\in M$ is called a groupless $F$-set based in $M$. We do allow in our definition of groupless $F$-sets $k=0$, in which case, the groupless $F$-set consists of the single point $b$. We denote by $\GF_M$ the set of all groupless $F$-sets based in $M$. For every subgroup $\Gamma\subset M$, we denote by $\GF_M(\Gamma)$ the collection of groupless $F$-sets contained in $\Gamma$ and based in $M$. When $M$ is clear from the context, we will drop the index $M$ from our notation.
\end{definition}

\begin{remark}
\label{R:base}
Each groupless $F$-set $O$ is based in a finitely generated $\mathbb{Z}[F]$-module.
\end{remark}

\begin{definition}
\label{D:F-sets}
Let $M$ be a module over a Frobenius ring $R[F]$. Let $\Gamma\subset M$ be a subgroup. A set of the form $(C+H)$, where $C\in\GF_M(\Gamma)$ and $H$ is a subgroup of $\Gamma$ invariant under $F$ is called an $F$-set in $\Gamma$ based in $M$. The collection of all such $F$-sets in $\Gamma$ is denoted by $\F_M(\Gamma)$. When $M$ is clear from the context, we will drop the index $M$ from our notation. 
\end{definition}

Let $G$ be a semiabelian variety defined over $\mathbb{F}_q$. Let $F$ be the corresponding Frobenius for $\mathbb{F}_q$. Let $K$ be a finitely generated regular extension of $\mathbb{F}_q$. We fix an algebraic closure $K^{\alg}$ of $K$. Let $\Gamma$ be a finitely generated subgroup of $G(K)$. We denote by $\F(\Gamma)$ and $\GF(\Gamma)$ the collection of $F$-sets and respectively, the collection of groupless $F$-sets in $\Gamma$ based in $G(K^{\alg})$ (which is obviously a $\mathbb{Z}[F]$-module). When we do not mention the $\mathbb{Z}[F]$-submodule containing the base points for the $F$-sets contained in $\Gamma$, then we will always understand that the corresponding submodule is $G(K^{\alg})$. The following theorem is our main result for semiabelian varieties.

\begin{theorem}
\label{T:groups}
Let $G$, $K$ and $\Gamma$ be as in the above paragraph. Let $X$ be a $K$-subvariety of $G$. Then $X(K)\cap\Gamma=\bigcup_{i=1}^r(C_i+H_i)$, where $(C_i+\Delta_i)\in\F(\Gamma)$. Moreover, the subgroups $\Delta_i$ from $X(K)\cap\Gamma$ are of the form $G_i(K)\cap\Gamma$, where $G_i$ is an algebraic subgroup of $G$ defined over $\mathbb{F}_q$.
\end{theorem}

As mentioned in Section~\ref{S:intro}, the result of our Theorem~\ref{T:groups} was establised in \cite{F-sets} (see Theorem $7.8$) and in \cite{revisited} (see Theorem $2.1$) for finitely generated $\mathbb{Z}[F]$-modules $\Gamma\subset G(K)$. Because $\mathbb{Z}[F]$ is a finite extension of $\mathbb{Z}$, each finitely generated $\mathbb{Z}[F]$-module is also a finitely generated group (but \emph{not} every finitely generated group is invariant under $F$).

We describe now the setting for our Drinfeld modules statements. We start by defining Drinfeld modules over finite fields.

Let $p$ be a prime number and let $q$ be a power of $p$. Let $C$ be a projective nonsingular curve defined over $\mathbb{F}_q$. We fix a closed point $\infty$ on $C$. Let $A$ be the ring of $\mathbb{F}_q$-valued functions on $C$, regular away from $\infty$. Then $A$ is a Dedekind domain. Moreover, $A$ is a finite extension of $\mathbb{F}_q[t]$. Hence, for every nonzero prime ideal $\pp\subset A$, $A/\pp$ is a finite field.

Let $F$ be the corresponding Frobenius on $\mathbb{F}_q$. 
We call a Drinfeld module defined over $\mathbb{F}_{q}$ a ring homomorphism $\phi:A\rightarrow\mathbb{F}_{q}[F]$ such that there exists $a\in A$ for which $\phi_a:=\phi(a)\notin\mathbb{F}_q\cdot F^0$ (i.e. the degree of $\phi_a$ as a polynomial in $F$ is positive). In general, for every $a\in A$, we write $\phi_a$ to denote $\phi(a)\in\mathbb{F}_q[F]$. We note that this is not the most general definition for Drinfeld modules defined over finite fields (see Example~\ref{E:sharp}).

For each field extension $L$ of $\mathbb{F}_{q}$, $\phi$ induces an action on $\mathbb{G}_a(L)$ by $a*x:=\phi_a(x)$ for every $x\in L$ and for every $a\in A$. For each $g\ge 1$, we extend the action of $A$ diagonally on $\mathbb{G}_a^g$.

Clearly, for every $a\in A$, $F\phi_a=\phi_aF$. This means $F$ is an \emph{endomorphism} of $\phi$ (see Section $4$ of Chapter $2$ in \cite{Gos}). We let $A[F]\in\End(\phi)$ be the finite extension of $A$ generated by $F$, where we identified $A$ with its image in $\mathbb{F}_q[F]$ through $\phi$. Actually, $A[F]$ is isomorphic to $\mathbb{F}_q[F]$. However, we keep the notation $A[F]$ instead of $\mathbb{F}_q[F]$, when we talk about modules over this ring only to emphasize the Drinfeld module action given by $A$.

\begin{lemma}
\label{L:Frobenius ring}
The ring $A[F]$ defined in the above paragraph is a Frobenius ring.
\end{lemma}

\begin{proof}
Because for some $a\in A$, $\phi_a$ is a polynomial in $F$ of positive degree, we conclude $F$ is integral over $A$. Because $\mathbb{F}_q[F]$ is a domain, we conclude $F$ is not a zero divisor. Also, no nonzero element of $A[F]$ is infinitely divisible by $F$ because all elements of $A[F]$ are polynomials in $F$ and so, no nonzero polynomial can be infinitely divisible by some polynomial of positive degree. Therefore $A[F]$ is a Frobenius ring.
\end{proof}

Let $K$ be a regular field extension of $\mathbb{F}_{q}$. We fix an algebraic closure $K^{\alg}$ of $K$. Let $\Gamma$ be a finitely generated $A[F]$-submodule of $\mathbb{G}_a^g(K)$. We denote by $\F(\Gamma)$ and $\GF(\Gamma)$ the $F$-sets and respectively, the groupless $F$-sets in $\Gamma$ based in $\mathbb{G}_a^g(K^{\alg})$. When we do not mention the $A[F]$-submodule containing the base points for the $F$-sets contained in $\Gamma$, we will always understand that the corresponding submodule is $\mathbb{G}_a^g(K^{\alg})$. We will explain in Section~\ref{S:Drinfeld} that the following Mordell-Lang statement for Drinfeld modules defined over finite fields follows along the same lines as Theorem $7.8$ in \cite{F-sets}.

\begin{theorem}
\label{T:DF-sets}
Let $\phi:A\rightarrow\mathbb{F}_q[F]$ be a Drinfeld module. Let $K$ be a regular extension of $\mathbb{F}_{q}$. Let $g$ be a positive integer. Let $\Gamma$ be a finitely generated $A[F]$-submodule of $\mathbb{G}_a^g(K)$ and let $X$ be an affine subvariety of $\mathbb{G}_a^g$ defined over $K$. Then $X(K)\cap\Gamma$ is a finite union of $F$-sets in $\Gamma$.
\end{theorem}

\section{The Mordell-Lang Theorem for semiabelian varieties defined over finite fields}
\label{S:semiabelian}

\begin{proof}[Proof of Theorem~\ref{T:groups}.] 
We first observe that the subgroups $\Delta_i$ from the intersection of $X$ with $\Gamma$ are indeed of the form $G_i(K)\cap\Gamma$ for algebraic groups $G_i$ defined over $\mathbb{F}_q$. Otherwise, we can always replace a subgroup $\Delta_i$ appearing in the intersection $X(K)\cap\Gamma$ with its Zariski closure $G_i$ and then intersect with $\Gamma$ (see also the proof of Lemma $7.4$ in \cite{F-sets}). Because $G_i$ is the Zariski closure of a subset of $G(K)$, then $G_i$ is defined over $K$. Because $G_i$ is an algebraic subgroup of $G$, then $G_i$ is defined over $\mathbb{F}_q^{\alg}$. Because $K$ is a regular extension of $\mathbb{F}_q$, we conclude that $G_i$ is defined over $\mathbb{F}_q=K\cap\mathbb{F}_q^{\alg}$.

We will prove the main statement of Theorem~\ref{T:groups} by induction on $\dim(X)$. Clearly, when $\dim(X)=0$ the statement holds (the intersection is a finite collection of points in that case). Assume the statement holds for $\dim(X)<d$ and we prove that it holds also for $\dim(X)=d$.

We will use in our proof a number of reduction steps.

\emph{Step 1.} Because $X(K)\cap\Gamma=\overline{X(K)\cap\Gamma}\cap\Gamma$ we may assume that $X(K)\cap\Gamma$ is Zariski dense in $X$.

\emph{Step 2.} At the expense of replacing $X$ by one of its irreducible components, we may assume $X$ is irreducible. Each irreducible component of $X$ has Zariski dense intersection with $\Gamma$. If our Theorem~\ref{T:groups} holds for each irreducible component of $X$, then it also holds for $X$.

\emph{Step 3.} We may assume the stabilizer $\Stab_G(X)$ of $X$ in $G$ is finite. Indeed, let $H:=\Stab_G(X)$. Then $H$ is defined over $K$ (because $X$ is defined over $K$) and also, $H$ is defined over $\mathbb{F}_q^{\alg}$ (because it is an algebraic subgroup of $G$). Thus $H$ is defined over $\mathbb{F}_q$. Let $\pi:G\rightarrow G/H$ be the natural projection. Let $\hat{G}$, $\hat{X}$ and $\hat{\Gamma}$ be the images of $G$, $X$ and $\Gamma$ through $\pi$. Clearly $\hat{\Gamma}$ is a finitely generated subgroup of $\hat{G}(K)$ and also, $\hat{X}$ is defined over $K$. 

If $\dim(H)>0$, then $\dim(\hat{X})<\dim(X)=d$. Hence, by the inductive hypothesis, $\hat{X}(K)\cap\hat{\Gamma}$ is a finite union of $F$-sets in $\hat{\Gamma}$. Using the fact that the kernel of $\pi|_{\Gamma}$ stabilizes $X(K)\cap\Gamma$, we conclude
$$X(K)\cap\Gamma=\pi|_{\Gamma}^{-1}\left(\hat{X}(K)\cap\hat{\Gamma}\right),$$
which shows that $X(K)\cap\Gamma$ is also a finite union of $F$-sets, because $\ker\left(\pi|_{\Gamma}\right)$ is a subgroup of $\Gamma$ invariant under $F$ (we recall that $\ker(\pi)=H$ is invariant under $F$). 

Therefore, we work from now on under the assumptions that 

$(i)$ $\overline{X(K)\cap\Gamma}=X$;

$(ii)$ $X$ is irreducible;

$(iii)$ $\Stab_G(X)$ is finite.

Let $\Gaa$ be the $\mathbb{Z}[F]$-module generated by $\Gamma$. Because $\Gamma$ is finitely generated and $F$ is integral over $\mathbb{Z}$, then also $\Gaa$ is finitely generated. By Theorem $7.8$ of \cite{F-sets}, $X(K)\cap\Gaa$ is a finite union of $F$-sets in $\Gaa$. So, there are finitely many groupless $F$-sets $C_i$ and $\mathbb{Z}[F]$-submodules $H_i\subset\Gaa$ such that 
$$X(K)\cap\Gaa=\cup_i\left(C_i+H_i\right).$$
We want to show $\bigcup_i\left(C_i+H_i\right)\cap\Gamma$ is a finite union of $F$-sets in $\Gamma$. It suffices to show that for each $i$, there exists a finite union $B_i$ of $F$-sets in $\Gamma$ such that $\left(C_i+H_i\right)\cap\Gamma\subset B_i\subset X(K)$. Indeed, the existence of such $B_i$ yields
$$X(K)\cap\Gamma=\cup_i B_i,$$ as desired.

\emph{Case 1.} $\dim\overline{C_i+H_i}<d$.

Let $X_i:=\overline{C_i+H_i}$. Then $X_i$ is defined over $K$ (because $(C_i+H_i)\subset G(K)$) and $\dim(X_i)<d$. So, by the induction hypothesis, $B_i:=X_i(K)\cap\Gamma$ is a finite union of $F$-sets in $\Gamma$. Clearly, $\left(C_i+H_i\right)\cap\Gamma\subset B_i\subset X(K)$ (because $X_i\subset X$).

\emph{Case 2.} $\dim\overline{C_i+H_i}=d$.

Because $X=\overline{X(K)\cap\Gamma}$, then $X=\overline{X(K)\cap\Gaa}$. Moreover, $X$ is irreducible and so, because $\dim\overline{C_i+H_i}=\dim (X)$, then $X=\overline{C_i+H_i}$. Hence $H_i\subset\Stab_G(X)$ because $$C_i+H_i+H_i=C_i+H_i\text{ and so, $\overline{C_i+H_i}+H_i\subset\overline{C_i+H_i}$.}$$ Because $\Stab_G(X)$ is finite, we conclude $H_i$ is finite. Thus $(C_i+H_i)$ is a finite union of groupless $F$-sets because it can be written as a finite union $\cup_{h\in H_i}\left(h+C_i\right)$. We let $B_i:=\left(C_i+H_i\right)\cap\Gamma$. We will show that for each (of the finitely many elements) $h\in H_i$, 
\begin{equation}
\label{E:cycle 1}
\left(h+C_i\right)\cap\Gamma\text{ is a finite union of groupless $F$-sets in $\Gamma$.}
\end{equation}
The following lemma will prove \eqref{E:cycle 1} and so, it will conclude the proof of Theorem~\ref{T:groups}.

\begin{lemma}
\label{L:orbit}
Let $M$ be a finitely generated $\mathbb{Z}[F]$-submodule of $G(K^{\alg})$ and let $O\in\GF_M$. If $\Gamma$ is a finitely generated subgroup of $G(K^{\alg})$, then $O\cap\Gamma$ is a finite union of groupless $F$-sets based in $M$.
\end{lemma}

\begin{proof}
If $O\cap\Gamma$ is finite, then we are done. So, from now on, we may assume $O\cap\Gamma$ is infinite. Also, we may and do assume $\Gamma\subset M$ (otherwise we replace $\Gamma$ with $\Gamma\cap M$).

Let $O:=Q+S(P_1,\dots,P_k;\delta_1,\dots,\delta_k)$, where $Q,P_1,\dots,P_k\in M$ and $\delta_i\in\mathbb{N}^{*}$ for every $i\in\{1,\dots,k\}$. We may assume that $\delta_1=\dots=\delta_k=1$, in which case $S(P_1,\dots,P_k;\delta_1,\dots,\delta_k):=S(P_1,\dots,P_k;1)$. Indeed, if we show that $$\left(Q+S(P_1,\dots,P_k;1)\right)\cap\Gamma\text{ is a union of groupless $F$-sets,}$$ then also its subsequent intersection with $\left(Q+S(P_1,\dots,P_k;\delta_1,\dots,\delta_k)\right)$ is a finite union of groupless $F$-sets, as shown in part $(a)$ of Lemma $3.7$ in \cite{F-sets}.

Because $M$ is a finitely generated $\mathbb{Z}$-module, $M$ is isomorphic with a direct sum of its finite torsion $M_{\tor}$ and a free $\mathbb{Z}$-submodule $M_1$. 

Let 
\begin{equation}
\label{E:minimal polynomial}
f(X):=X^g-\sum_{i=0}^{g-1}\alpha_i X^i
\end{equation}
be the minimal polynomial for $F$ over $\mathbb{Z}$ (i.e. $f(F)=0$ in $\End(G)$). Let $r_1,\dots,r_g$ be all the roots in $\mathbb{C}$ of $f(X)$. 
Clearly, each $r_i\ne 0$ because $F$ is not a zero-divisor in $\End(G)$. Also, each $r_i$ has absolute value larger than $1$ (actually, their absolute values equal $q$ or $q^{\frac{1}{2}}$, according to the Riemann hypothesis for semiabelian varieties defined over $\mathbb{F}_q$). Finally, all $r_i$ are distinct. At most one of the $r_i$ is real and it equals $q$ (and it corresponds to the multiplicative part of $G$), while all of the other $r_i$ have absolute value equal to $q^{\frac{1}{2}}$ (and they correspond to the abelian part of $G$). If $$0\rightarrow T\rightarrow G\rightarrow A\rightarrow 0$$ is a short exact sequence of group varieties, with $T$ being a torus and $A$ an abelian variety, then the roots $r_i$ of absolute value $q^{\frac{1}{2}}$ correspond to roots of the minimal polynomial over $\mathbb{Z}$ for the Frobenius morphism on $A$. The abelian variety $A$ is isogenuous with a product of simple abelian varieties $A_i$, all defined over a finite field. If $f_i$ is the minimal polynomial of the corresponding Frobenius on $A_i$, then the minimal polynomial $f_0$ of the Frobenius on $A$ is the least common multiple of all $f_i$. For each $i$, $\End(A_i)$ is a domain and so, $f_i$ has simple roots. Therefore $f_0$ (and so, $f$) has simple roots.

The definition of $f$ shows that for every point $P\in G(K^{\alg})$, 
\begin{equation}
\label{E:recursive}
F^gP=\sum_{j=0}^{g-1}\alpha_jF^jP.
\end{equation}
We conclude that for all $n\ge g$,
\begin{equation}
\label{E:recursive 1}
F^nP=\sum_{j=0}^{g-1}\alpha_jF^{n-g+j}P.
\end{equation}

For each $j$ we define the sequence $\{z_{j,n}\}_{n\ge 0}$ as follows
\begin{equation}
\label{E:1}
z_{j,n}=0\text{ if $0\le n\le g-1$ and $n\ne j$};
\end{equation}
\begin{equation}
\label{E:2}
z_{j,j}=1\text{ and}
\end{equation}
\begin{equation}
\label{E:3}
z_{j,n}=\sum_{l=0}^{g-1}\alpha_lz_{j,n-g+l}\text{ for all $n\ge g$.}
\end{equation}
Using \eqref{E:1} and \eqref{E:2} we obtain that
\begin{equation}
\label{E:4}
F^nP=\sum_{j=0}^{g-1}z_{j,n}F^jP\text{, for every $0\le n\le g-1$.}
\end{equation}
We prove by induction on $n$ that 
\begin{equation}
\label{E:relation}
F^nP=\sum_{j=0}^{g-1}z_{j,n}F^jP\text{, for every $n\ge 0$.}
\end{equation}
We already know \eqref{E:relation} is valid for all $n\le g-1$ due to \eqref{E:4}. Thus we assume \eqref{E:relation} holds for all $n< N$, where $N\ge g$ and we prove that \eqref{E:relation} also holds for $n=N$.
Using \eqref{E:recursive 1}, we get
\begin{equation}
\label{E:6}
F^NP=\sum_{j=0}^{g-1}\alpha_jF^{N-g+j}.
\end{equation}
We apply the induction hypothesis to all $F^{N-g+j}$ for $0\le j\le g-1$ and conclude
\begin{equation}
\label{E:7}
\sum_{j=0}^{g-1}\alpha_jF^{N-g+j}=\sum_{j=0}^{g-1}\alpha_j\sum_{i=0}^{g-1}z_{i,N-g+j}F^iP=\sum_{i=0}^{g-1}\left(\sum_{j=0}^{g-1}\alpha_jz_{i,N-g+j}\right)F^iP.
\end{equation}
We use \eqref{E:3} in \eqref{E:7} and conclude
\begin{equation}
\label{E:8}
\sum_{j=0}^{g-1}\alpha_jF^{N-g+j}=\sum_{i=0}^{g-1}z_{i,N}F^iP.
\end{equation}
Combining \eqref{E:6} and \eqref{E:8} we obtain the statement of \eqref{E:relation} for $n=N$. This concludes the inductive proof of \eqref{E:relation}.

Because $\{z_{j,n}\}_n$ is a recursive defined sequence, then for each $j\in\{0,\dots,g-1\}$ there exist
$\{\gamma_{j,l}\}_{1\le l\le g}\subset\mathbb{Q}^{\alg}$ such that for every $n\in\mathbb{N}$,
\begin{equation}
\label{E:gammas}
z_{j,n}=\sum_{1\le l\le g}\gamma_{j,l}r_l^n.
\end{equation}
To derive \eqref{E:gammas} we also use the fact that all $r_i$ are distinct, nonzero numbers.

Equations \eqref{E:relation} and \eqref{E:gammas} show that for every $n$ and for every $P\in G(K^{\alg})$,
\begin{equation}
\label{E:needed}
F^nP=\sum_{0\le j\le g-1}\left(\sum_{1\le l\le g}\gamma_{j,l}r_l^n\right)F^jP.
\end{equation}
 
For each $i\in\{1,\dots,k\}$ and for each $j\in\{0,\dots,g-1\}$, let $F^jP_i:=T^{(j)}_i+Q^{(j)}_i$, with $T^{(j)}_i\in M_{\tor}$ and $Q^{(j)}_i\in M_1$. Also, let $Q:=T_0+Q_0$, where $T_0\in M_{\tor}$ and $Q_0\in M_1$. 

Let $R_1,\dots,R_m$ be a basis for the $\mathbb{Z}$-module $M_1$. For each $j\in\{0,\dots,g-1\}$ and for each $i\in\{1,\dots,k\}$, let 
\begin{equation}
\label{E:coefficient}
Q^{(j)}_i:=\sum_{l=1}^m a_{i,j}^{(l)}R_l,
\end{equation}
where $a_{i,j}^{(l)}\in\mathbb{Z}$.
Finally, let $a_0^{(1)},\dots,a_{0}^{(m)}\in\mathbb{Z}$ such that $Q_0=\sum_{j=1}^m a_{0}^{(j)}R_j$. 

For every $n\in\mathbb{N}$ and for every $i\in\{1,\dots,k\}$, \eqref{E:relation} and the definitions of $Q^{(j)}_i$ and $T^{(j)}_i$ yield
\begin{equation}
\label{E:the sum}
F^nP_i=\sum_{0\le j\le g-1}z_{j,n}\left(T^{(j)}_i+Q^{(j)}_i\right)=\sum_{0\le j\le g-1}z_{j,n} T^{(j)}_i+\sum_{0\le j\le g-1}z_{j,n} Q^{(j)}_i.
\end{equation}
Because $T^{(j)}_i\in M_{\tor}$, then for each $(n_1,\dots,n_k)\in\mathbb{N}^k$, $$\sum_{i=1}^k \sum_{j=0}^{g-1} z_{j,n_i}T^{(j)}_i\in M_{\tor}.$$ Also, because $Q_0$ and all $Q^{(j)}_i$ are in $M_1$ and because $z_{j,n}\in\mathbb{Z}$, then for each $(n_1,\dots,n_k)\in\mathbb{N}^k$, 
$$Q_0+\sum_{i=1}^k\sum_{j=0}^{g-1} z_{j,n_i}Q^{(j)}_i\in M_1.$$
Moreover, 
\begin{equation}
\label{E:point splitting}
Q+\sum_{i=1}^kF^{n_i}P_i=\left(T_0+\sum_{\substack{1\le i\le k\\ 0\le j\le g-1}}z_{j,n_i}T^{(j)}_i\right)+\left(Q_0+\sum_{\substack{1\le i\le k\\ 0\le j\le g-1}}z_{j,n_i}Q^{(j)}_i\right).
\end{equation}

For each $h\in M_{\tor}$, if $(h+M_1)\cap\Gamma$ is not empty, we fix $(h+U_h)\in\Gamma$ for some $U_h\in M_1$. Let $\Gamma_1:=\Gamma\cap M_1$. Then 
\begin{equation}
\label{E:h}
\left(h+M_1\right)\cap\Gamma=h+U_h+\Gamma_1.
\end{equation}
For each $h\in M_{\tor}$, we let $O_h:=\{P\in O\mid P=h+P'\text{ with $P'\in M_1$}\}$. Then using \eqref{E:h}, we get
\begin{equation}
\label{E:small intersection}
O\cap\Gamma=\bigcup_{h\in M_{\tor}}O_h\cap\left(h+U_h+\Gamma_1\right)=\bigcup_{h\in M_{\tor}}\left(h+\left((-h+O_h)\cap(U_h+\Gamma_1)\right)\right).
\end{equation}
Clearly, $(-h+O_h)\in M_1$. Therefore \eqref{E:small intersection} and \eqref{E:point splitting} yield
\begin{equation}
\label{E:big intersection}
O\cap\Gamma=\bigcup_{\substack{h\in M_{\tor}\\T_0+\sum_{i,j}z_{j,n_i}T^{(j)}_i=h}}\left(h+\left(\left(Q_0+\sum_{i,j}z_{j,n_i}Q^{(j)}_i\right)\bigcap\left(U_h+\Gamma_1\right)\right)\right).
\end{equation}
In \eqref{E:big intersection}, the union is over the finitely many torsion points of $M_{\tor}$ ($M$ is finitely generated) and it might be that not for each $h\in M_{\tor}$ there is a corresponding nonempty intersection in \eqref{E:big intersection}.

Fix $h\in M_{\tor}$. We show that the set of tuples $(n_1,\dots,n_k)\in\mathbb{N}^k$ for which
\begin{equation}
\label{E:h_0}
h=T_0+\sum_{i,j}z_{j,n_i}T^{(j)}_i
\end{equation}
is a finite union of cosets of \emph{semigroups} of $\mathbb{N}^k$ (a semigroup of $\mathbb{N}^k$ is the intersection of a subgroup of $\mathbb{Z}^k$ with $\mathbb{N}^k$). Indeed, let $N\in\mathbb{N}^{*}$ such that $M_{\tor}\subset G[N]$. Because for each $j\in\{0,\dots,g-1\}$, $z_{j,n}$ is a recursively defined sequence (as shown by \eqref{E:1}, \eqref{E:2} and \eqref{E:3}), then the sequence $\{z_{j,n}\}_n$ is eventually periodic modulo $N$ (a recursively defined sequence is eventually periodic modulo any integral modulus). Thus each value taken by $T_0+\sum_{i,j}z_{j,n_i}T^{(j)}_i$ is attained for tuples $(n_1,\dots,n_k)$ which belong to a finite union of cosets of semigroups of $\mathbb{N}^k$.

We will prove next that for each fixed $h\in M_{\tor}$, the tuples $(n_1,\dots,n_k)$ for which 
\begin{equation}
\label{E:nontorsion}
\left(Q_0+\sum_{i,j}z_{j,n_i}Q^{(j)}_i\right)\in\left(U_h+\Gamma_1\right)
\end{equation} 
form a finite union of cosets of semigroups of $\mathbb{N}^k$. This will finish the proof of our theorem because this result, combined with the one from the previous paragraph and combined with \eqref{E:big intersection}, will show that the tuples $(n_1,\dots,n_k)$ for which
$$Q+\sum_{i=1}^k F^{n_i}P_i\in\Gamma$$
form a finite union of cosets of semigroups of $\mathbb{N}^k$ (we are also using the fact that the intersection of two finite unions of cosets of semigroups is also a finite union of cosets of semigroups). Lemma $3.4$ of \cite{F-sets} shows that the set of points in $ O$ corresponding to a finite union of cosets of semigroups containing the tuples of exponents $(n_1,\dots,n_k)$ is a finite union of groupless $F$-sets.

Because $\Gamma_1\subset M_1$ and $M_1$ is a free $\mathbb{Z}$-module with basis $\{R_1,\dots,R_m\}$, we can find (after a possible relabelling of $R_1,\dots,R_m$) a $\mathbb{Z}$-basis $V_1,\dots,V_n$ ($n\le m$) of $\Gamma_1$ of the following form:
$$V_1=\beta_{1}^{(i_1)}R_{i_1}+\dots+\beta_{1}^{(m)}R_m;$$
$$V_2=\beta_{2}^{(i_2)}R_{i_2}+\dots+\beta_{2}^{(m)}R_m;$$
and in general, $V_j=\beta_{j}^{(i_j)}R_{i_j}+\dots+\beta_{j}^{(m)}R_m$, where $$1\le i_1<i_2<\dots<i_n\le m$$ and all $\beta_{j}^{(i)}\in\mathbb{Z}$. We also assume $\beta_j^{(i_j)}\ne 0$ for every $j\in\{1,\dots,n\}$ ($n\ge 1$ because we assumed the intersection $O\cap\Gamma$ is infinite, which means $\Gamma_1$ is infinite, because otherwise $|O\cap\Gamma |\le |M_{\tor}|$).

Let $b_{0}^{(1)},\dots,b_{0}^{(m)}\in\mathbb{Z}$ such that $U_h=\sum_{j=1}^m b_0^{(j)}R_j$. Then a point $$P:=\sum_{j=1}^m c^{(j)}R_j\in (U_h+\Gamma_1)$$ if and only if there exist integers $k_1,\dots,k_n$ such that 
\begin{equation}
\label{E:v1}
P=U_h+\sum_{i=1}^nk_iV_i.
\end{equation}
Using the expressions of the $V_i$, $U_h$ and $P$ in terms of the $\mathbb{Z}$-basis $\{R_1,\dots,R_m\}$ of $M_1$, we obtain the following relations:
\begin{equation}
\label{E:a}
c^{(j)}=b_0^{(j)}\text{ for every $1\le j<i_1$;}
\end{equation}
\begin{equation}
\label{E:v2}
c^{(j)}=b_0^{(j)}+k_1\beta_1^{(j)}\text{ for every $i_1\le j<i_2$;}
\end{equation}
\begin{equation}
\label{E:v3}
c^{(j)}=b_0^{(j)}+k_1\beta_1^{(j)}+k_2\beta_2^{(j)}\text{ for every $i_2\le j<i_3$}
\end{equation}
and so on, until
\begin{equation}
\label{E:v4}
c^{(m)}=b_0^{(m)}+\sum_{i=1}^nk_i\beta_i^{(m)}.
\end{equation}
We express equation \eqref{E:v2} for $j=i_1$ as a linear congruence modulo $\beta_1^{(i_1)}$ and obtain
\begin{equation}
\label{E:b}
c^{(i_1)}\equiv b_0^{(i_1)}\left(\text{ mod $\beta_1^{(i_1)}$}\right).
\end{equation}
Also from \eqref{E:v2} for $j=i_1$, we get $k_1=\frac{c^{(i_1)}-b_0^{(i_1)}}{\beta_1^{(i_1)}}$. Then we substitute this formula for $k_1$ in \eqref{E:v2} for all $i_1<j<i_2$ and obtain
\begin{equation}
\label{E:c}
c^{(j)}=b_0^{(j)}+\frac{c^{(i_1)}-b_0^{(i_1)}}{\beta_1^{(i_1)}}\beta_1^{(j)}\text{ for every $i_1<j<i_2$.}
\end{equation}
Then we express \eqref{E:v3} for $j=i_2$ as a linear congruence modulo $\beta_2^{(i_2)}$ (also using the expression for $k_1$ computed above). We obtain
\begin{equation}
\label{E:d}
c^{(i_2)}\equiv b_0^{(i_2)}+\frac{c^{(i_1)}-b_0^{(i_1)}}{\beta_1^{(i_1)}}\beta_1^{(i_2)}\left(\text{ mod $\beta_2^{(i_2)}$}\right).
\end{equation}
Next we equate $k_2$ from \eqref{E:v3} for $j=i_2$ (also using the formula for $k_1$) and obtain $$k_2=\frac{c^{(i_2)}-b_0^{(i_2)}-\frac{c^{(i_1)}-b_0^{(i_1)}}{\beta_1^{(i_1)}}\beta_1^{(i_2)}}{\beta_2^{(i_2)}}$$ Then we substitute this formula for $k_2$ in \eqref{E:v3} for $i_2<j<i_3$ and obtain
\begin{equation}
\label{E:z}
c^{(j)}=b_0^{(j)}+\frac{c^{(i_1)}-b_0^{(i_1)}}{\beta_1^{(i_1)}}\cdot\beta_1^{(j)}+\frac{c^{(i_2)}-b_0^{(i_2)}-\frac{c^{(i_1)}-b_0^{(i_1)}}{\beta_1^{(i_1)}}\beta_1^{(i_2)}}{\beta_2^{(i_2)}}\cdot\beta_2^{(j)}.
\end{equation}
We go on as above until we express $c^{(m)}$ in terms of $$c^{(i_1)},\dots,c^{(i_n)}$$ and $b_0^{(m)}$ and the $\beta_j^{(l)}$.
We observe that all congruences can be written as linear congruences over $\mathbb{Z}$. For example, the above congruence equation \eqref{E:d} modulo $\beta_2^{(i_2)}$ can be written as the following linear congruence over $\mathbb{Z}$:
$$\beta_1^{(i_1)}\cdot c^{(i_2)}\equiv \left(c^{(i_1)}-b_0^{(i_1)}\right)\beta_1^{(i_2)}+\beta_1^{(i_1)}b_0^{(i_2)}\left(\text{ mod $\beta_1^{(i_1)}\cdot\beta_2^{(i_2)}$}\right).$$
Hence all the above conditions are either linear congruences or linear equations for the $c^{(j)}$. 

A typical intersection point from the inner intersection in \eqref{E:big intersection} corresponding to a tuple $(n_1,\dots,n_k)\in\mathbb{N}^k$ is $$\left(Q_0+\sum_{i,j} z_{j,n_i}Q^{(j)}_i\right)\cap\left(U_h+\Gamma_1\right)$$ and it can be written in the following form (see also \eqref{E:gammas}):
$$\sum_{l=1}^g \left(a_0^{(l)}+\sum_{\substack{1\le i\le k\\0\le j\le g-1}}a_{i,j}^{(l)}\sum_{e=1}^{g}\gamma_{j,e}r_{e}^{n_i}\right)R_l.$$
Such a point lies in $(U_h+\Gamma_1)$ if and only if its coefficients 
$$a_0^{(l)}+\sum_{\substack{1\le i\le k\\0\le j\le g-1}}a_{i,j}^{(l)}\sum_{e=1}^{g}\gamma_{j,e}r_{e}^{n_i}$$
with respect to the $\mathbb{Z}$-basis $\{R_1,\dots,R_m\}$ of $M_1$ satisfy the linear congruences and linear equations such as \eqref{E:a}, \eqref{E:b}, \eqref{E:c}, \eqref{E:d} and \eqref{E:z}, associated to $(U_h+\Gamma_1)$. A linear equation as above yields an equation of the following form (after collecting the coefficients of $r_e^{n_i}$ for each $1\le e\le g$ and each $1\le i\le k$):
\begin{equation}
\label{E:M-L}
\sum_{e=1}^g\sum_{i=1}^k d_{e,i} r_e^{n_i}=D.
\end{equation}
All $d_{e,i}$ and $D$ are algebraic numbers. A tuple $(n_1,\dots,n_k)\in\mathbb{N}^k$ satisfying \eqref{E:M-L} corresponds to an intersection point of the linear variety $L$ in $\left(\mathbb{G}_m^{g}\right)^k(\mathbb{Q}^{\alg})$ given by the equation
\begin{equation}
\label{E:line}
\sum_{e=1}^g\sum_{i=1}^k d_{e,i} X_{e,i}=D
\end{equation}
and the finitely generated subgroup $G_0$ of $\left(\mathbb{G}_m^{g}\right)^k(\mathbb{Q}^{\alg})$ spanned by 
\begin{equation}
\label{E:vector}
(r_1,\dots,r_g,1,\dots,1)\text{; }(1,\dots,1,r_1,\dots,r_g,1,\dots,1)\text{; }\dots, (1,\dots,1,r_1,\dots,r_g).
\end{equation}
Each vector in \eqref{E:vector} has $gk$ components.
There are $k$ multiplicatively independent generators above for $G_0$ (we are using the fact that $|r_i|>1$, for each $i$). Hence $G_0\isom\mathbb{Z}^k$. By Lang Theorem for $\mathbb{G}_m^{gk}$, we conclude the intersection of $L(\mathbb{Q}^{\alg})$ and $G_0$ is a finite union of cosets of subgroups of $G_0$. The subgroups of $G_0$ correspond to subgroups of $\mathbb{Z}^k$. Hence the tuples $(n_1,\dots,n_k)\in\mathbb{N}^k$ which satisfy \eqref{E:M-L} belong to a finite union of cosets of semigroups of $\mathbb{N}^k$.

A congruence equation as \eqref{E:b} or \eqref{E:d}, corresponding to conditions for a point to lie in $(U_h+\Gamma_1)$ yields a congruence relation between the coefficients (with respect to the $\mathbb{Z}$-basis $\{R_1,\dots,R_m\}$ of $M_1$)
of a typical point of the form $Q_0+\sum_{j,i}z_{j,n_i}Q^{(j)}_i$. We will show that such tuples $(n_1,\dots,n_k)$ belong to a finite union of cosets of semigroups of $\mathbb{N}^k$.

The coefficient of $R_l$ in $\left(Q_0+\sum_{j,i}z_{j,n_i}Q^{(j)}_i\right)$ can be written as (see also \eqref{E:coefficient})
\begin{equation}
\label{E:coefficient 2}
a_0^{(l)}+\sum_{\substack{1\le i\le k\\0\le j\le g-1}}a_{i,j}^{(l)}z_{j,n_i}.
\end{equation}
Hence a congruence equation corresponding to a point of the form $\left(Q_0+\sum_{j,i}z_{j,n}Q^{(j)}_i\right)$ which also lies in $(U_h+\Gamma_1)$ has the form
\begin{equation}
\label{E:congruence}
\sum_{j=0}^{g-1}\sum_{i=1}^k d_{j,i}z_{j,n_i}\equiv D_1(\text{ mod $D_2$})
\end{equation}
for some integers $d_{j,i}$ (we recall that $a_{i,j}^{(l)}\in\mathbb{Z}$), $D_1$ and $D_2\ne 0$.
Recursively defined sequences as $\{z_{j,n}\}_n$ are eventually periodic modulo any nonzero integer (hence, they are eventually periodic modulo $D_2$). Therefore all the solutions $(n_1,\dots,n_k)$ to \eqref{E:congruence} belong to a finite union of cosets of semigroups of $\mathbb{N}^k$.

Hence for each $h\in M_{\tor}$ the tuples $(n_1,\dots,n_k)\in\mathbb{N}^k$ for which $$\left(Q_0+\sum_{i,j}z_{j,n_i}Q^{(j)}_i\right)\in\left(U_h+\Gamma_1\right),$$
form a finite union of cosets of semigroups of $\mathbb{N}^k$. We also proved that for each $h\in M_{\tor}$ the tuples $(n_1,\dots,n_k)\in\mathbb{N}^k$ for which
$$h=T_0+\sum_{i,j}z_{j,n_i}T^{(j)}_i,$$
form a finite union of cosets of semigroups of $\mathbb{N}^k$.
In conclusion, we get that $$\left(Q+\sum_{i=1}^kF^{n_i}P_i\right)\in\Gamma$$ if and only if $(n_1,\dots,n_k)$ belongs to a finite union of cosets of semigroups of $\mathbb{N}^k$. The corresponding subset of $\left(Q+S(P_1,\dots,P_k;1)\right)$ for a finite union of cosets of semigroups of $\mathbb{N}^k$ is precisely a finite union of groupless $F$-sets based in $M$ (as shown by Lemma $3.4$ of \cite{F-sets}).

This concludes the proof of Lemma~\ref{L:orbit}
\end{proof}
As remarked before the statement of Lemma~\ref{L:orbit}, this lemma concludes the proof of our Theorem~\ref{T:groups}.
\end{proof}

\section{The Mordell-Lang Theorem for Drinfeld modules defined over finite fields}
\label{S:Drinfeld}

The setting for this section is that $\phi:A\rightarrow\mathbb{F}_q[F]$ is a Drinfeld module.

The following result (which is the equivalent for Drinfeld modules of Lemma $7.5$ in \cite{F-sets}) will be used in the proof of our Theorem~\ref{T:DF-sets}.
\begin{lemma}
\label{L:F-cosets}
Let $K$ be a finitely generated field extension of $\mathbb{F}_q$ and let $\Gamma\subset\mathbb{G}_a^g(K)$ be a finitely generated $A[F]$-submodule. 

$(a)$ The $F$-pure hull of $\Gamma$ in $\mathbb{G}_a^g(K)$, i.e. the set of all $x\in\mathbb{G}_a^g(K)$ such that $F^mx\in\Gamma$ for some $m\ge 0$, is a finitely generated $A$-module. In particular, $\Gamma$ is a finitely generated $A$-module.

$(b)$ For each $m>0$, $\Gamma/F^m\Gamma$ is finite.

$(c)$ There exists $m\ge 0$ such that $\Gamma\setminus F\Gamma\subset\mathbb{G}_a^g(K)\setminus\mathbb{G}_a^g\left(K^{q^{m}}\right)$.
\end{lemma}

\begin{proof}
$(a)$ First we observe that the $F$-pure hull $\Gaa$ of $\Gamma$ is an $A[F]$-module, and so, implicitly an $A$-module. Indeed, if $x\in\Gaa$ and $m\in\mathbb{N}$ such that $F^mx\in\Gamma$, then for every $f\in A[F]$, $$F^m(f(x))=f(F^mx)\in f(\Gamma)\subset\Gamma.$$
Therefore $f(x)\in\Gaa$, showing that $\Gaa$ is an $A[F]$-module.

It suffices to prove $(a)$ under the extra assumption that $\Gamma=\Gamma_0^g$ (the cartesian product of $\Gamma_0$ with itself $g$ times), where $\Gamma_0\subset K$ is a finitely generated $A[F]$-module. Indeed, let $\Gamma_0$ be the finitely generated $A[F]$-submodule of $K$ spanned by all the generators (over $A[F]$) of the projections of $\Gamma$ on the $g$ coordinates of $\mathbb{G}_a^g(K)$. Clearly $\Gamma\subset\Gamma_0^g$ and if we prove $(a)$ for $\Gamma_0$, then the result of $(a)$ holds also for $\Gamma_0^g$ and implicitly for its submodule $\Gamma$ (the $F$-pure hull of $\Gamma$ is an $A$-submodule of the $F$-pure hull of $\Gamma_0^g$ and a submodule of a finitely generated module over a Dedekind domain is also finitely generated). So, we are left to show that the $F$-pure hull $\Gaa_0$ of $\Gamma_0$ in $K$ is a finitely generated $A$-module.

By its construction, $\Gamma_0$ is a finitely generated $A[F]$-submodule of $K$. Because $F$ is integral over $A$, we conclude $\Gamma_0$ is also finitely generated as an $A$-module. As explained in the beginning of our proof, $\Gaa_0$ is also an $A$-module. We first prove $\Gaa_0$ lies inside the $A$-division hull $\Gamma_0'$ of $\Gamma_0$ in $K$. Indeed, let $x\in\Gaa _0$ and let $m\in\mathbb{N}$ such that $F^mx\in\Gamma_0$. We will prove next that $x\in\Gamma_0'$.

Because $F$ is integral over $A$, then also $F^m$ is integral over $A$. Let $s\in\mathbb{N}^{*}$ and let $\alpha_0,\dots,\alpha_{s-1}\in A$ such that 
\begin{equation}
\label{E:s-minimal}
F^{sm}=\sum_{i=0}^{s-1}\alpha_iF^{i\cdot m}\text{ in $\End(\phi)$.}
\end{equation}
Because $A[F]$ is a domain, we may assume $\alpha_0\ne 0$ (otherwise we would divide \eqref{E:s-minimal} by powers of $F^m$ until the coefficient of $F^0$ would be nonzero). Equality \eqref{E:s-minimal} shows that 
\begin{equation}
\label{E:s-minimal partial}
\phi_{\alpha_0}(x)=F^{sm}x-\sum_{i=1}^{s-1}\phi_{\alpha_i}\left(F^{i\cdot m}x\right)\in\Gamma_0,
\end{equation}
because $F^mx\in\Gamma_0$ and $\Gamma_0$ is an $A[F]$-module. Thus \eqref{E:s-minimal partial} shows $x$ belongs to the $A$-division hull $\Gamma_0'$. Let $F_0:=\Frac(A)$. Because $\Gaa _0\subset\Gamma_0'$ and because $\Gamma_0$ is a finitely generated $A$-module, we conclude 
\begin{equation}
\label{E:rank inequality}
\dim_{F_0}\left(\Gaa _0\tensor_A F_0\right)\le\dim_{F_0}\left(\Gamma_0'\tensor_A F_0\right)<\aleph_0.
\end{equation}
Hence \eqref{E:rank inequality} shows $\Gaa _0$ has \emph{finite rank} as an $A$-module.
Lemma $4$ of \cite{poo} shows that every finite rank $A$-module is finitely generated. This concludes the proof of $(a)$.

$(b)$ Because $\Gamma$ is a finitely generated $A[F]$-module, then $\Gamma/F^m\Gamma$ is a finitely generated $A[F]/(F^m)$-module. Hence, it suffices to show $A[F]/(F^m)$ is a finite ring. Let, as before, \eqref{E:s-minimal} be the minimal equation of $F^m$ over $A$. Then $\alpha_0\in F^m\cdot A[F]$. So, $A[F]/(F^m)$ is a quotient of $A[F]/(\alpha_0)$. Clearly, $A[F]/(\alpha_0)\isom\left(A/(\alpha_0)\right)[F]$. Because $\alpha_0\ne 0$ and $A$ is a Dedekind domain for which the residue field for each nonzero  ideal is finite, we conclude $A/(\alpha_0)$ is finite (we know that $A/\pp$ is finite for every nonzero prime ideal $\pp$, but every nonzero ideal in $A$ is a product of nonzero prime ideals). Because $F$ is integral over $A$ we conclude $\left(A/(\alpha_0)\right)[F]$ is finite. Hence $A[F]/(F^m)$ is finite and so, $\Gamma/F^m\Gamma$ is finite, as desired.

$(c)$ Because the $F$-pure hull $\Gaa$ of $\Gamma$ in $\mathbb{G}_a^g(K)$ is finitely generated as an $A[F]$-module, then there exists $m_0>0$ such that $F^{m_0}\Gaa\subset\Gamma$. Let $m:=m_0+1$. Then $$\Gamma\cap\mathbb{G}_a^g\left(K^{q^{m}}\right)\subset F^m\Gaa\subset F\Gamma.$$ Hence $\Gamma\setminus F\Gamma\subset\mathbb{G}_a^g(K)\setminus\mathbb{G}_a^g\left(K^{q^{m}}\right)$.
\end{proof}

We will also use in our proof of Theorem~\ref{T:DF-sets} the following result on the combinatorics of the $F$-sets.
\begin{lemma}
\label{L:combinatorics 1}
Suppose $K$ is a regular field extension of $\mathbb{F}_q$, $\Gamma\subset\mathbb{G}_a^g(K)$ is a finitely generated $A[F]$-module, $X\subset\mathbb{G}_a^g$ is an affine variety defined over $K$ and $b\in\mathbb{N}^{*}$. Clearly $\Gamma$ is an $A[F^b]$-module as well. If $U\subset\Gamma$ is an $F^b$-set with $U\subset X(K)$, then there exists $V\in\F(\Gamma)$ such that $U\subset V\subset X(K)$. In particular, if $X(K)\cap\Gamma$ is a finite union of $F^b$-sets, then it is also a finite union of $F$-sets.
\end{lemma}

\begin{proof}
Our proof follows the proof of its similar statement for semiabelian varieties instead of Drinfeld modules and for $\mathbb{Z}[F]$ instead of $A[F]$ (Lemma $7.4$ of \cite{F-sets}).

Let $U=C+\Delta$, where $C$ is a groupless $F^b$-set and $\Delta$ is a subgroup of $\Gamma$ invariant under $F^b$. Let $H$ be the Zariski closure of $\Delta$ in $\mathbb{G}_a^g$. Then $H$ is invariant under $F^b$. Hence $H$ is defined over $\mathbb{F}_{q^b}$ (which is the fixed field of $F^b$). Because $H$ is the Zariski closure of a subset of $\mathbb{G}_a^g(K)$, then $H$ is defined over $K$. Therefore $H$ is defined over $K\cap\mathbb{F}_{q^b}$. Because $K$ is a regular extension of $\mathbb{F}_q$, then $K\cap\mathbb{F}_{q^b}=\mathbb{F}_q$. Thus $H$ is defined over $ \mathbb{F}_q$ and so, $H(K)\cap\Gamma$ is invariant under $F$.

Clearly every groupless $F^b$-set is also a groupless $F$-set and so, $C$ is a groupless $F$-set. Therefore we conclude that $V:=C+H(K)\cap\Gamma$ is an $F$-set in $\Gamma$, which contains $U$. On the other hand, $H\subset X$ (because $\Delta\subset X(K)$ and $H=\overline{\Delta}$). Moreover, for each $c\in C$, $$c+H(K)\cap\Gamma\subset\overline{c+\Delta}(K)\subset X(K).$$ Thus $V\subset X(K)$, as desired.
\end{proof}

The proof of the next two lemmas are identical with the proofs of Corollary $7.3$ and respectively, Lemma $3.9$ in \cite{F-sets}.
\begin{lemma}
\label{L:combinatorics 2}
Suppose $\Gamma\subset\mathbb{G}_a^g(K)$ is a finitely generated $A[F]$-module, $U$ is a finite union of $F$-sets in $\Gamma$ and $X\subset\mathbb{G}_a^g$ is an affine variety defined over $K$. Let $\Sigma:=\bigcup_{n\ge 0}F^nU$ and suppose that $\Sigma\subset X(K)$. Then there exists a finite union $B$ of $F$-sets in $\Gamma$ such that $\Sigma\subset B\subset X(K)$.
\end{lemma}

\begin{lemma}
\label{L:combinatorics 3}
Suppose $M$ is a finitely generated $A[F]$-module.

$(a)$ The intersection of two finite unions of $F$-sets in $M$ is also a finite union of $F$-sets in $M$.

$(b)$ If $X$ is a finite union of $F$-sets in $M$ and $N$ is a submodule of $M$, then $X\cap N$ is a finite union of $F$-sets in $N$.
\end{lemma}

We will deduce Theorem~\ref{T:DF-sets} from the following slightly more general statement (our Theorem~\ref{T:DF-sets} is a particular case of Theorem~\ref{T:DF-sets 2} for $H=\{0\}$).
\begin{theorem}
\label{T:DF-sets 2}
Let $K$ be a regular extension of $\mathbb{F}_q$. Let $H$ be any algebraic subgroup of $\mathbb{G}_a^g$ defined over $\mathbb{F}_q$. Then for every variety $X\subset\mathbb{G}_a^g/H$ defined over $K$ and for every finitely generated $A[F]$-submodule $\Gamma\subset\left(\mathbb{G}_a^g/H\right)(K)$, the intersection $X(K)\cap\Gamma$ is a finite union of $F$-sets in $\Gamma$ based in $\left(\mathbb{G}_a^g/H\right)(K^{\alg})$.
\end{theorem}

\begin{proof}
We first observe that because $H$ is an algebraic group defined over $\mathbb{F}_q$, then $H$ is invariant under $A[F]$. Hence, the quotient $\mathbb{G}_a^g/H$ is equipped with a natural $A$-action.

Our proof follows the proof of Theorem $7.8$ of \cite{F-sets}. Because $\phi$ is defined over a finite field and because $\Gamma$ is a finitely generated $A$-module (see $(a)$ of Lemma~\ref{L:F-cosets}) and because $X$ is defined over a finitely generated field, then there exists a finitely generated subfield $L$ of $K$ such that $X$ is defined over $L$ and $\Gamma\subset\mathbb{G}_a^g(L)$. Therefore we may and do assume that $K$ is finitely generated.

We will use induction on $\dim(X)$. If $\dim (X)=0$, then $X(K)\cap\Gamma$ is a finite collection of points. Clearly, each point is an $F$-set. We assume that Theorem~\ref{T:DF-sets 2} holds for $\dim(X)<n$ (for some $n\ge 1$) and we will prove that it also holds for $\dim(X)=n$.

We may assume $\overline{X(K)\cap\Gamma}=X$ (otherwise, we may replace $X$ with $\overline{X(K)\cap\Gamma}$). Also, we may assume $X$ is irreducible because it suffices to prove Theorem~\ref{T:DF-sets 2} for each irreducible component of $X$ (we are using the fact that the intersection of $X$ with $\Gamma$ is Zariski dense if and only if the intersection of each irreducible component of $X$ with $\Gamma$ is Zariski dense in that component).

The next lemma shows that a translate of $X$ is defined over a finite field. The proof of Lemma~\ref{L:over finite fields} is almost identical with the proof of Lemma $7.7$ in \cite{F-sets}. Lemma $7.7$ in \cite{F-sets} holds for any finitely generated subgroup of a semiabelian variety. In particular, it holds for any finitely generated $\mathbb{Z}[F]$-submodule of a semiabelian variety. The only difference between Lemma $7.7$ in \cite{F-sets} and our Lemma~\ref{L:over finite fields} is that in \cite{F-sets}, $\Gamma$ can be taken to be a module over the Frobenius ring $\mathbb{Z}[F]$ (associated to a semiabelian variety defined over a finite field), while in our case, $\Gamma$ is a module over the Frobenius ring $A[F]$ (associated to a Drinfeld module defined over a finite field). The only property of the Frobenius ring used in the proof of Lemma $7.7$ in \cite{F-sets} is property $(b)$ from Lemma~\ref{L:F-cosets} and the only property of the ambient algebraic group $G$ (a semiabelian variety in \cite{F-sets} and $\mathbb{G}_a^g/H$ for us) used in the proof of Lemma $7.7$ in \cite{F-sets} is that $\bigcap_{n\ge 1}F^n G(K^{\alg})=G(\mathbb{F}_q^{\alg})$.
\begin{lemma}
\label{L:over finite fields}
Suppose $\Gamma$ is a finitely generated $A[F]$-submodule of $\left(\mathbb{G}_a^g/H\right)(K)$ and $X\subset\left(\mathbb{G}_a^g/H\right)$ is a variety defined over $K$ such that $X(K)\cap\Gamma$ is Zariski dense in $X$. Then for some $\gamma\in K^{\alg}$, $(\gamma+X)$ is defined over $\mathbb{F}_q^{\alg}$.
\end{lemma}

Next we show that we may assume $X$ is defined over $\mathbb{F}_q$. Lemma~\ref{L:over finite fields} shows that there exists $\gamma\in K^{\alg}$ such that $(\gamma+X)$ is defined over $\mathbb{F}_q^{\alg}$. Let $\Gamma'$ be the finitely generated $A[F]$-module generated by $\gamma$ and the elements of $\Gamma$. Let $K':=K(\gamma)$. Because $X(K)\cap\Gamma$ is Zariski dense in $X$, then $(\gamma+X)\cap\Gamma'$ is Zariski dense in $(\gamma+X)$. Hence $(\gamma+X)$ is defined over $K'$. But we already know that $(\gamma+X)$ is defined over $\mathbb{F}_q^{\alg}$. Hence $(\gamma+X)$ is defined over $$\mathbb{F}_{q^b}:=K'\cap\mathbb{F}_q^{\alg}.$$
Assuming the statement of our Theorem~\ref{T:DF-sets 2} valid for varieties defined over the finite field fixed by the Frobenius, we obtain that $(\gamma+X)\cap\Gamma'$ is an $F^b$-set. Because $\Gamma$ is an $A[F^b]$-submodule of $\Gamma'$, we conclude $$X(K)\cap\Gamma=X(K')\cap\Gamma=\left(X(K')\cap\Gamma'\right)\cap\Gamma.$$
Hence, using part $(b)$ of Lemma~\ref{L:combinatorics 3}, $X(K)\cap\Gamma$ is an $F^b$-set in $\Gamma$. An application of Lemma~\ref{L:combinatorics 1} concludes the proof that $X(K)\cap\Gamma$ is indeed an $F$-set in $\Gamma$. Therefore, from now on, we assume that $X$ is defined over $\mathbb{F}_q$.

We may also assume $\Stab(X)\subset\mathbb{G}_a^g/H$ is trivial. Indeed, let $H_1=\Stab(X)$. Then $H_1$ is defined over the same field as $X$. Hence $H_1$ is defined over $\mathbb{F}_q$. We consider the canonical quotient map $\pi:\left(\mathbb{G}_a^g/H\right)\rightarrow\mathbb{G}_a^g/(H+H_1)$. Let $\hat{X}$ and $\hat{\Gamma}$ be the images of $X$ and $\Gamma$ through $\pi$. Clearly $\Stab(\hat{X})=\{0\}$. Moreover, if Theorem~\ref{T:DF-sets 2} holds for $\hat{X}(K)\cap\hat{\Gamma}$, then it also holds for $X(K)\cap\Gamma=\pi|_{\Gamma}^{-1}\left(\hat{X}(K)\cap\hat{\Gamma}\right)$ (we use the fact that $\ker\left(\pi|_{\Gamma}\right)=\Gamma\cap H_1(K)$ is a subgroup of $\Gamma$ invariant under $F$). Also, it is precisely this part of our proof where we need the hypothesis of Theorem~\ref{T:DF-sets 2} be that $X$ is a subvariety of a quotient of $\mathbb{G}_a^g$ through an algebraic subgroup defined over $\mathbb{F}_q$.

From this point on the proof of Theorem~\ref{T:DF-sets 2} is identical with the proof of Theorem $7.8$ in \cite{F-sets} (we provided in Lemmas~\ref{L:F-cosets}, \ref{L:combinatorics 1} and \ref{L:combinatorics 2} the technical ingredients that are used in the argument from the proof of Theorem $7.8$ in \cite{F-sets}).
\end{proof}

The following result follows from Theorem $3.1$ in \cite{revisited} the same way our Theorem~\ref{T:DF-sets} followed from Theorem $7.8$ in \cite{F-sets}.
\begin{theorem}
\label{T:revisited}
Let $\phi:A\rightarrow\mathbb{F}_q[F]$ be a Drinfeld module. Let $F$ be the Frobenius on $\mathbb{F}_q$. Let $K$ be an algebraically closed field extension of $\mathbb{F}_q$. Let $X\subset\mathbb{G}_a^g$ (for some $g\ge 1$) be an affine variety defined over $K$. Let $\Gamma\subset\mathbb{G}_a^g(K)$ be a finitely generated $A[F]$-module. Let $\Gamma':=\Gamma+\mathbb{G}_a^g(\mathbb{F}_q^{\alg})$. Then $X(K)\cap\Gamma'$ is a finite union of sets of the form $\left(U+Y(\mathbb{F}_q^{\alg})\right)$, where $U\subset\Gamma'$ is an $F^b$-set for some $b\in\mathbb{N}^{*}$ and $Y\subset\mathbb{G}_a^g$ is an affine variety defined over $\mathbb{F}_q^{\alg}$. 
\end{theorem}

In the following Example~\ref{E:sharp}, we extend the notion of Drinfeld modules defined over finite fields and then we show that for our \emph{new} Drinfeld modules, the groups appearing in the intersection from the conclusion of Theorem~\ref{T:DF-sets} are not necessarily $A$-modules (and hence, they are not $A[F]$-modules). This is in contrast with the semiabelian case where the groups appearing in the intersection $X(K)\cap\Gamma$ are $\mathbb{Z}[F]$-modules.

\begin{example}
\label{E:sharp}
Let $a\in\mathbb{N}^{*}$. Let $K$ be a regular extension of $\mathbb{F}_{q^a}$. Let $\mathbb{F}_{q^a}\{F\}$ be the ring of twisted polynomials in $F$ with coefficients in $\mathbb{F}_{q^a}$ (the addition is the usual one, while the multiplication is the composition of functions). A Drinfeld module over a finite field is a ring homomorphism $\phi:A\rightarrow\mathbb{F}_{q^a}\{F\}$ for which there exists $a\in A$ such that $\phi_a\notin\mathbb{F}_{q^a}\cdot F^0$. Then $F$ is not necessarily an endomorphism for $\phi$, but $F^a\in\End(\phi)$. We want to characterize the intersections $X(K)\cap\Gamma$, where $X\subset\mathbb{G}_a^g$ is an affine variety defined over $K$ and $\Gamma\subset\mathbb{G}_a^g(K)$ is a finitely generated $A[F^a]$-submodule.

We cannot always expect that the subgroups of $\Gamma$ appearing in $X(K)\cap\Gamma$ be actually $A$-submodules. For example, let $C=\mathbb{P}_{\mathbb{F}_q}^1$ and let $A=\mathbb{F}_q[t]$. Let $a=2$. Define $\phi:A\rightarrow\mathbb{F}_{q^2}\{F\}$ by $\phi_t=F+F^3$. Let $\lambda\in\mathbb{F}_{q^2}\setminus\mathbb{F}_q$. Consider the curve $X\subset\mathbb{G}_a^2$ defined by $y=\lambda x$. Let $K=\mathbb{F}_{q^2}(t)$ and let $\Gamma\subset\mathbb{G}_a^2(K)$ be the cyclic $A[F^2]$-submodule generated by $(t,\lambda t)$. 

Then $X(K)\cap\Gamma$ consists of all points in $\Gamma$ of the form $\left(f(t),f(\lambda t)\right)$, where $f\in A[F^2]$ is a polynomial in $F^2$ with coefficients in $\mathbb{F}_{q}$. In particular, $X(K)\cap\Gamma$ is invariant under $\phi_{t^2}=F^2+2F^4+F^6$, but it is not invariant under $\phi_t$. So, the intersection is an $\mathbb{F}_q[F^2]$-submodule of $\Gamma$, but not an $A[F^2]$-submodule.
\end{example}

The following example shows that we cannot obtain a similar statement as our Theorem~\ref{T:groups} in the context of Drinfeld modules, i.e. we cannot replace the $A[F]$-submodules $\Gamma$ in Theorem~\ref{T:DF-sets} with simply $A$-modules.
\begin{example}
\label{E:finite Drinfeld}
Assume $q$ is odd and let $A=\mathbb{F}_q[t]$. Define $\phi:A\rightarrow\mathbb{F}_q[F]$ by $\phi_t=F+F^2$.

Let $Y\subset\mathbb{G}_a^g$ be a smooth curve defined over $\mathbb{F}_q$ and let $K:=\mathbb{F}_q(Y)$. Let $P\in Y(K)$ be a generic point for $Y$. Define $X:=\overline{Y+Y}$ and assume $X$ does not contain translates of nontrivial algebraic subgroups of $\mathbb{G}_a^g$ (for \emph{generic} curves $Y$ this is always possible). Let $\Gaa$ be the $A[F]$-submodule of $\mathbb{G}_a^g(K)$ generated by $P$. Then, using that $X$ does not contain a translate of a nontrivial algebraic subgroup of $\mathbb{G}_a^g$, we conclude
\begin{equation}
\label{E:contain}
X(K)\cap\Gaa=S(P,P;1).
\end{equation}
Let $\Gamma$ be the cyclic $A$-module generated by $P$. Clearly, $\Gamma\subset\Gaa$. Hence, using \eqref{E:contain}, we obtain
$$X(K)\cap\Gamma=\bigcup_{n\ge 0}\phi_{t^{p^n}}(P)=\bigcup_{n\ge 0}\left(F^{p^n}P+F^{2p^n}P\right).$$ This is the case because the only elements $a\in A$ such that $\phi_a=F^n+F^m$ are of the form $a=t^{p^n}$ (this is an easy exercise in combinatorics, whose proof we provide below for completeness).
\end{example}

\begin{lemma}
\label{L:example}
Assume $p$ is an odd prime and let $q$ be a power of $p$. Let $A:=\mathbb{F}_q[t]$ and define the Drinfeld module $\phi:A\rightarrow\mathbb{F}_q[F]$ by $\phi_t=F+F^2$. Then the only elements $a\in A$ such that $\phi_a$ equals $F^n+F^m$ for some $n,m\in\mathbb{N}$ are of the form $a=t^{p^n}$ (in which case $\phi_{t^{p^n}}=F^{p^n}+F^{2p^n}$.
\end{lemma}

\begin{proof}
Let $a=\sum_{i=0}^n a_it^i\in A$ (hence $a_i\in\mathbb{F}_q$). Assume $\phi_a$ is the sum of two powers of $F$. We will prove that all $a_i=0$ for $i<n$ and also that $n$ is a power of $p$.

First we observe that if $a_i=0$ for all $i<n$, then $a=a_n t^n$ and so, the expansion of $(F+F^2)^n$ contains \emph{only} two powers of $F$ if and only if $n$ is a power of $p$ (Lucas Theorem for Binomial Congruences). Moreover, $a_n=1$ in order for $\phi_a$ to be a sum of two powers of $F$.

Assume there is $k<n$ such that $a_k\ne 0$. Let $m$ be the least such $k$. Then the term $a_m F^m$ has the smallest power of $F$ which appears in $\phi_a$ (and it is not cancelled by any other term in $\phi_a$). On the other hand, $a_nF^{2n}$ is the term in $\phi_a$ with the largest power of $F$ (and also it is not cancelled by any other term in $\phi_a$). Therefore the only two powers of $F$ in $\phi_a$ are $F^m$ and $F^{2n}$.

Let $l$ be the index of the first nonzero digit in the expansion of $n$ in base $p$, i.e. $$n=\sum_{j\ge l}\alpha_jp^j$$ and $\alpha_l\ne 0$. Then the coefficient of $F^{2n-p^l}$ in the expansion $\phi_{a_nt^n}$ equals $a_n\binom{n}{p^l}\ne 0$ in $\mathbb{F}_q$ (by Lucas Theorem for Binomial Congruences). Moreover, also by Lucas Theorem, we get that $F^{2n-p^l}$ is the largest power of $F$, not equal to $F^{2n}$, which appears with nonzero coefficient in the expansion of $\phi_{a_nt^n}$. Also, $$2n-p^l\ge n>m.$$ Thus the power $F^{2n-p^l}$ has to be cancelled by another term in $\phi_a$. Let $n_1<n$ be the largest index $i$ such that $a_i\ne 0$. Then the largest power of $F$ in $\phi_{a-a_nt^n}$ is $F^{2n_1}$ which does not cancel $F^{2n-p^l}$, because $p^l$ is odd. Hence, either the power $F^{2n-p^l}$ or the power $F^{2n_1}$ appear with nonzero coefficients in $\phi_a$, contradicting thus the fact that the only powers of $F$ in $\phi_a$ are $F^m$ and $F^{2n}$.
\end{proof}

\begin{remark}
The above proof works applied to the Drinfeld module $\phi:\mathbb{F}_q[t]\rightarrow\mathbb{F}_q[F]$ defined by $\phi_t=F+F^3$, in case $p=2$, and shows that the only elements $a\in A$ such that $\phi_a$ equals $F^n+F^m$ for some $n,m\in\mathbb{N}$ are of the form $a=t^{2^n}$ (in which case $\phi_{t^{2^n}}=F^{2^n}+F^{3\cdot 2^n}$). This allows us to construct a similar example in characteristic $2$ as Example~\ref{E:finite Drinfeld} for the failure of a Mordell-Lang statement such as Theorem~\ref{T:DF-sets} for finitely generated $A$-modules $\Gamma$.
\end{remark}

\end{document}